\documentclass[12pt,a4paper]{amsart}
\usepackage{amsmath,amssymb,epsf,graphics}
\usepackage{amsfonts,amsxtra,amscd,verbatim,eucal}
\usepackage{enumerate}
\usepackage[all]{xy}

\newtheorem{thm}{Theorem}[section]
\newtheorem{lemma}{Lemma}
\newtheorem{prop}{Proposition}

\newtheorem{defn}{Definition}

\theoremstyle{remark}
\newtheorem{rem}{Remark}

\def\R{{\mathbb R}}

\begin{document}
\author{Catriona Maclean\\
Institut Fourier\\ Universit\'e Grenoble Alpes.}
\title{Characterising approximable algebras.}
\maketitle
\begin{abstract}
In \cite{volumes_chen}, Huayi Chen introduces the notion of an approximable 
graded algebra, which he uses to prove a Fujita-type theorem in the arithmetic 
setting,  and asked if any such algebra is the graded ring of a big line 
bundle on a projective variety. This was proved to be false in \cite{maclean}. 
Continuing the analysis
started in \cite{maclean}, we here show that whilst not every approximable
graded algebra is a sub algebra of the graded ring of a big line bundle 
on a projective variety, it is the case that to any 
approximable graded algebra ${\bf B}$ we can associate a projective divisor $X({\bf B})$
and an infinite divisor $D({\bf B}) =\sum_{i=1}^\infty a_i D_i$ with $a_i\rightarrow 0$ 
such that ${\bf B}$ is included in 
\[ R( D({\bf B}))=\oplus_n H^0(X({\bf B}), n D({\bf B})).\]
We also establish a partial converse to these results by showing that if the
infinite divisor $D=\sum_i a_iD_i$ converges in the space of numerical classes 
then any full-dimensional sub-graded algebra of 
$\oplus_mH^0(X, \lfloor mD \rfloor))$ is approximable.
\end{abstract}
\section{Introduction}
The Fujita approximation theorem, \cite{fuj}, is an important result in 
algebraic geometry. It states that whilst the section ring associated to a 
big line bundle $L$ on an algebraic variety $X$ 
\[ R(L)\stackrel{\rm def}{=} \oplus_m H^0(mL, X)\]
is typically not a finitely generated algebra, it can be approximated
arbitrarily well by finitely generated algebras. More precisely,  we have that
\begin{thm}[Fujita]
Let $X$ be an algebraic variety and let $L$ be a big line bundle on $X$. 
For any $\epsilon>0$ there exists a birational modification 
\[\pi: \hat{X}\rightarrow X\]
and a decomposition of $\mathbb{Q}$ divisors, $\pi^* (L)= A+E$ such that
\begin{itemize}
\item $A$ is ample and $E$ is effective,
\item ${\rm vol}(A)\geq (1-\epsilon){\rm vol}(L)$.
\end{itemize}
\end{thm}
In \cite{LM}, Lazarsfeld and Mustata used the Newton-Okounkov body associated to
$A$ to give a simple proof of Fujita approximation. The Newton-Okounkov body, 
constructed in \cite{KK} and \cite{LM}, building on previous work of Okounkov
\cite{Okounkov}, is a convex body $\Delta_{Y_\bullet}(L,X)$ in 
$\mathbb{R}^d$ associated to the 
data of
\begin{itemize}
\item a $d$-dimensional variety $X$
\item an admissible flag $Y_\bullet$ on $X$
\item a big line bundle $L$ on $X$.
\end{itemize}
This convex body encodes information
on the asymptotic behaviour of the spaces of global sections $H^0(nL)$ for large
values of $L$. \\ \\
Lazarsfeld and 
Mustata's simple proof of Fujita approximation is based on the equality of 
volumes of Newton-Okounkov bodies
\begin{equation}\label{volumes} {\rm vol}(L)= d! {\rm vol}
(\Delta_{Y_\bullet}(L,X))\end{equation}
where we recall that the volume of a big line bundle on a $d$-dimesional variety
is defined by
\[ {\rm vol}(L)= \lim_{n\rightarrow \infty}\frac{ d! h^0(nL)}{n^d}.\]
One advantage of their approach to the Fujita theorem
is that Newton-Okounkov bodies are not 
only defined for section algebras $R(L)$, but also for any graded sub-algebra of
section algebras. Lazarsfeld and Mustata give combinatorical 
conditions (conditions 2.3-2.5 of \cite{LM}) under which equation
\ref{volumes} holds for a graded sub-algebra ${\bf B}= \oplus_m B_m\subset R(L)$
and show that these conditions hold if 
the graded subalgebra ${\bf B}$ 
{\it contains an ample series}.\footnote{Ie. 
if there exists an ample divisor $A\leq L$ such that $\oplus_m H^0(\lfloor
mA \rfloor)\subset B$}\\ \\
Di Biagio and Pacenzia in \cite{dBP}
subsequently used Newton-Okounkov bodies associated to 
restricted algebras to prove a Fujita approximation theorem for restricted 
linear series, ie. subalgebras of $\oplus_m H^0(mL|_V,V)$ obtained as the
restriction of the complete algebra $\oplus_m H^0(mL,X)$, where $V\subset X$ is a 
subvariety. \\ \\
In \cite{volumes_chen}, Huayi Chen uses Lazarsfeld and Mustata's work on 
Fujita approximation to
prove a Fujita-type approximation theorem in the arithmetic setting.
In the course of this work he defines the notion of approximable 
graded algebras, which are exactly those algebras for which
a Fujita-type approximation theorem hold. 
\begin{defn}
An integral graded algebra ${\bf B}=\oplus_m B_m$ with $B_0=k$ a field 
is approximable if and only if
the following conditions are satisfied.
\begin{enumerate}
\item all the graded pieces $B_m$ are finite dimensional over $k$.
\item for all sufficiently large $m$ the space $B_m$ is non-empty
\item  
for any $\epsilon$ there exists an $p_0$ such that for all $p\geq p_0$ 
we have that 
\[ \liminf_{n\rightarrow \infty} \frac{ {\rm dim}({\rm Im}(S^nB_p\rightarrow B_{np}))}{{\rm dim}(B_{np})}> (1-\epsilon). \]
\end{enumerate}
\end{defn}
In his paper \cite{volumes_chen} 
Chen asks 
whether any graded approximable algebra is in fact a subalgebra of the algebra 
of sections of a big line bundle. 
A counter-example was given to this is \cite{maclean}, where a counter example is constructed in which the 
graded approximable algebra is equal to the section ring of an {\it infinite} 
divisor\footnote{Infinite in this context meaning an infinite sum of Weil 
divisors with real coefficients
$\sum_i a_i D_i$.}. This begs the question : is
any approximable algebra a subalgebra of the section ring of an infinite divisior ?\\ \\
In the current paper we will prove that the answer is yes by establishing the following theorem.
\begin{thm}
Let ${\bf B}=\oplus_m B_m$ be a graded approximable algebra whose first graded piece $B_0$ is an algebraically closed field of characteristic zero. 
There is then a projective variety $X({\bf B})$ and an infinite divisor $D({\bf B})=\sum_{i=1}^\infty a_i D_i$ such that $a_i\rightarrow 0$ and there
is a natural inclusion of graded algebras
\[ {\bf B}\hookrightarrow \oplus_m H^0(X({\bf B}), m D({\bf B})).\]
\end{thm}
Furthermore, in the other direction we prove the following.
\begin{thm}
Let $X$ be a complex algebraic variety and let $D=\sum a_i D_i$ be an infinite 
Weil divisor on $X$ such that the sum of divisor classes 
$ \sum_i a_i[D_i ]$ converges
to a finite real big cohomology class. Any graded subalgebra of 
$\oplus_m H^0(mD)$ such that 
\[ \left( \frac{{\rm rk}(B_m)}{ m^{d({\bf C})}}\right)\]
does not converge to zero is then an approximable algebra.
\end{thm}
\section{Notation and two preliminary results.}
\noindent
In this section we will fix some notation and recall an essential preliminary lemma from \cite{volumes_chen}.
\\ \\
Throughout this article, $k$ will be an algebraically closed field of charcateristic zero. ${\bf B}=\oplus_m B_m$
will be a graded approximable algebra such that $B_0=k$: we will say that 
${\bf B}$ is a graded approximable algebra over $k$.  
For any natural numbers $k$ and $n$ we will 
denote by ${\rm Sym}^n(B_k)$ the $n$-th symmetric power of the vector space $B_k$ and by $S^n(B_k)$
the image of ${\rm Sym}^n(B_k)$ in $B_{nk}$. \\ \\
For any $k$ we will denote by $\langle B_k\rangle$ the subalgebra $\oplus_n S^n(B_k)\subset {\bf B}$ and
by ${\bf B}_k$ the subalgebra $\oplus_n B_{nk}$. \\ \\
We now recall a result from Chen on approximable algebras which will be 
necessary in what follows.
\begin{lemma}\label{volexist}[Chen, Proposition 2.4]
Let ${\bf B}=\oplus_{m\geq 0} B_m$ be an integral graded algebra which is
approximable. There then exists a constant $a\in \mathbb{N}^*$ such that, for 
any sufficiently large integer $p$, the algebra $\langle B_p\rangle$ has Krull 
dimension $a$. Furthermore, let us denote by 
denote by $d({\bf B})$ the number $a-1$. The sequence $v_n$ defined by
\begin{equation}\label{eqlim}v_n=\left( \frac{{\rm rk} B_n}{ n^{d({\bf B})}/ d({\bf B})}\right)_{n\geq 1}\end{equation}
then converges in $\mathbb{R}_+$
\end{lemma}
Naturally, the number $d({\bf B})$ represent the dimension of the algebra ${\bf B}$
and the limit of the sequence in equation \ref{eqlim} is its volume.
\begin{defn}
Let ${\bf B}=\oplus_{m\geq 0} B_m$ be an integral graded algebra over a field $k$ which is
approximable. We then define the {\bf dimension} of ${\bf B}$ to be the 
number $d({\bf B})$ whose existence is guaranteed by Lemma \ref{volexist}. Furthermore, we define the {\bf volume } of ${\bf B}$, denoted ${\rm vol}({\bf B})$ 
to be the limit 
\[{\rm vol}({\bf B})= \lim_{n\rightarrow \infty}
\left( \frac{{\rm rk} B_n}{ n^{d({\bf B})}/ d({\bf B})}\right)\]
\end{defn}
We will say that a graded, not a priori approximable, algebra 
$\oplus_m B_m$ {\it is of dimension $d$ and has volume $v$} if the sequence
\[ \lim_{n\rightarrow \infty}
\left( \frac{{\rm rk} B_n}{ n^{d({\bf B})}/ d({\bf B})}\right)\]
converges to the real number $v$. \\ \\
Note that the condition $(3)$ in the definition of an approximable alegbra tells
us that 
\[ \lim_{p\rightarrow \infty}({\rm vol}(\langle B_p\rangle)= 
{\rm vol}({\bf B}).\]
The following lemma from Chen (Corollary 2.5) will also be useful.
\begin{lemma}\label{chenlim}
Let ${\bf B}$ be a graded approximable algebra. We then have that for any $r\in \mathbb{N}$
\[ \lim_{n\rightarrow\infty}\frac{{\rm rk}(B_{n+r})}{ {\rm rk}(B_n)}=1.\]
\end{lemma}
\section{Construction of the variety $X({\bf B})$ and the divisors $D_m$.}
In the current section we will construct the variety $X({\bf B})$ and divisors $D_m$ associated to the
linear series $B_m$.
\subsection{Construction of $X({\bf B})$}
We will define the variety $X$ using the homogeneous field of fractions of ${\bf B}=\oplus_m B_m$,
which we now define.
\begin{defn}
Let ${\bf B}=\oplus_m B_m$ be a graded algebra such that $B_0=k$. Then we define its homogeneous 
fraction field by
\[ K^{\rm hom}({\bf B})= \left\{ \frac{b_1}{b_2}| \exists m \mbox{ such that } b_1, b_2\in B_m, b_2\neq 0\right\}/ \sim \]
where $\sim$ is the equivalence relation 
\[ \frac{b_1}{b_2}\sim \frac{c_1}{c_2} \Leftrightarrow b_1 c_2= c_1 b_2.\]
\end{defn}
\noindent
Note that $k$ is included in $K^{\rm hom}({\bf B})$ via the map 
$\lambda\rightarrow \frac{\lambda f}{f}$ for any $f\in B_m$. \\ \\
Choose $m$ large enough that $B_n$ and $B_{n+1}$ are both non-trivial. Choose
$f_1\in B_n$ and $f_2\in B_{n+1}$. For any $m$ we can then identify $B_m$ with a 
subspace of $K^{\rm hom}({\bf B})$ via the identification
\[ b_m\rightarrow \frac{b_m f_1^m}{f_2^m}.\]
Thoughout what follows, we will consider the space $B_m$ as a subvector space in $K^{\rm hom}({\bf B})$. \\ \\ 
The idea of our definition is that $X({\bf B})$ 
will be an algebraic variety whose 
function field will be $K^{\rm hom}({\bf B})$. 
Before being able to pose this definition, we will need to show that 
$K^{\rm hom}({\bf B})$ is finitely generated as a field 
extension of $k$.  
\begin{prop}\label{propfg}
Let ${\bf B}$ be an approximable graded algebra over an algebraically closed 
field $k$ of characteristic zero. The field $K^{\rm hom}({\bf B})$ is then a 
finitely generated field over $k$, whose
transcendence degree is equal to the dimension $d({\bf B})$. 
\end{prop}
\noindent
{\bf Proof of Proposition \ref{propfg}.}\\ \\
Suppose that ${\bf B}$ is an approximable graded algebra, and let $p_0$ be such that for any $p> p_0$ we have that 
\[ \liminf_n \left(\frac{{\rm dim}(S^n B_p)}{{\rm dim}(B_{np})}\right)>\frac{2}{3}\] 
We claim that $K^{\rm hom}({\bf B })$ is then generated as a field by $B_p\subset K^{\rm hom}({\bf B})$.\\ \\
Indeed, consider an arbitrary element 
\[\frac{b_1}{b_2} \in K^{\rm hom}({\bf B})\] with $b_1, b_2 \in B_m$
for some suficiently large $m$. After multiplication by an element of the form
$f/f$
we may assume that $m=kp$ is divisible by $p$. For $n$ large enough we have that 
\[ {\rm dim} \left(b_1\cdot {S}^{n}(B_p)\right) > 
\left(\frac{2}{3} {\rm dim}(B_{np})\right)> \left(\frac{2}{3}-\epsilon \right) {\rm dim}(B_{np+m})\]
where the last equality follows from Lemma \ref{chenlim}.
Similarly, we have that 
\[   {\rm dim} \left({\rm Sym}^{n+k}(B_p) \right)> \left(\frac{2}{3} {\rm dim}(B_{np+m})\right)\]
and hence the space \[S^{n+k}(B_p)\cap  b_1\cdot S^n(B_p)\] is of 
strictly positive dimension. Take
a non-zero element $b_3$ of this space: we then have that
\[ b_3= b_1 P_1 = P_2\]
for some elements $P_1\in S^n(B_p)$ and $P_2\in S^{n+k}(B_p)$. Similarly, there 
are elements
$Q_1\in S^{n}(B_p)$ and $Q_2\in S^{n+k}(B_p)$ such that  $b_2Q_1= Q_2$. But it
then follows that 
\[ \frac{b_1}{b_2}= \frac{P_2Q_1}{P_1Q_2}\]
in $K^{\rm hom}({\bf B})$, and since $P_1, P_2, Q_1,Q_2$
are all generated by $B_p$ this completes the proof of the first part of
Proposition 
\ref{propfg}. \\ \\
Indeed, the above proof establishes not only that 
$K^{\rm hom}({\bf B})$ is finitely 
generated, but moreover that it is equal to $K^{\rm  hom}(\langle B_p\rangle)$ 
for any sufficiently large $p$. \\ \\
It remains to show that the transcendence degree of $K^{\rm hom}({\bf B})$ 
is equal to the dimension of ${\bf B}$ as an approximable algebra. 
Since the algebra $\langle B_p \rangle$ is a finitely generated algebra over 
$k$, taking $p$ large enough if necessary  we have by \cite{Eisenbud}, 
Theorem A (p. 223), that 
\[ {\rm trdeg}_kK\left(\langle B_p \rangle\right)= a\]
Note that the above field is the total field of fractions of 
$\langle B_p\rangle$, not just the homogeneous part. We have that 
\[ K(\langle B_p \rangle)=K^{\hom}(\langle B_p \rangle)(f)\]
for any function $f\in K^{\rm hom}(\langle B_p \rangle)$ of the form
\[ f= \frac{f_2}{f_1}\]
with $f_1\in B_n$ and $f_2\in B_{n+1}$.  Moreover, $f$ is transcendent over 
$K^{\hom}(\langle B_p\rangle)$ by degree considerations, so it follows that
\[ {\rm trdeg}_kK^{\rm hom}\left(\langle B_p \rangle \right)= a-1=d({\bf B}).\]
This completes the proof of Proposition \ref{propfg}. \\ \\
Now, we are ready to define the variety $X({\bf B})$.
\begin{defn}
The variety $X({\bf B})$ is a smooth projective k-variety such that the function field
$K\left(X({\bf B})\right)= K^{\rm hom}({\bf B})$.
\end{defn}
\begin{rem}
The variety $X({\bf B})$ is here defined only up to birational equivalence. 
It can of course be chosen 
smooth by Hironaka's resolution of singularities.
\end{rem}
\begin{rem}
Since the dimension of any algebraic variety $X$ is equal to the transcendence 
degree of its function field, we have that 
\[ d({\bf B})= {\rm dim}(X).\]
\end{rem}
\begin{rem}
It follows from the proof of Proposition \ref{propfg} that 
$K\left(X({\bf B})\right)= 
K(\langle B_p \rangle )$, and in particular, the map defined on 
$X({\bf B})$ by the linear series $B_p$
is birational onto its image.
\end{rem} 
\begin{defn}
For any $b_m\in B_m$ we denote by $(b_m)_X$ the principal divisor on 
$X({\rm B})$ cut out by the rational function $b_m$. It negative part will be
denoted by $(b_m)_X^-$ and its positive part will be denoted by $(b_m)_X^+$, 
so that
\[ (b_m)_X= (b_m)_X^+-(b_m)_X^-.\]
\end{defn}
\subsection{Construction of $D_m$.}
The divisor $D({\bf B})$ will be constructed as the limit of the sequence of 
divisors 
$D_m/m$, where the divisors $D_m$ are constructed as poles of the rational functions 
$b_m\in B_m$.
\begin{defn}
For any $m$ such that $B_m$ is non-empty 
we define the effective divisor $D_m$ on $X({\bf B})$ by
\[ D_m= \sup_{b_m\in B_m}\left((b_m)_X^-\right).\]
where the supremum is taken with respect to the natural partial order on ${\rm Weil}(X({\bf B}))$.
\end{defn}
\noindent
We note that for any $b_m, b_m'\in B_m$ and for generic $\lambda\in k$ we 
have
\[ (b_m+\lambda b'_m)_X^-= \sup\left( (b_m)^-_X, (b'_m)^-_X\right)\]
so this supremum is actually a maximum. It follows that $D_m$ is indeed a finite divisor. \\ \\
Another possible characterisation of $D_m$ is the following
\begin{defn}
$D_m$ is the smallest divisor on $X({\bf B})$ with the property that $D_m+(b_m)
\geq 0$ for all $b_m\in B_m$.
\end{defn}
The construction of $D$ as the limit of the normalised divisors $D_m/m$ will 
depend on the following lemma.
\begin{lemma}\label{mono}
Let $m_1$ and $m_2$ be two natural numbers, and let the divisors $D_m$ be
as constructed above. 
If $m_1|m_2$ then $\frac{D_{m_1}}{m_1}\leq \frac{D_{m_2}}{m_2}$
\end{lemma}
\noindent
{\bf Proof of Lemma \ref{mono}}\\ \\
We have that $D_{m_1}= (b_{m_1})^-_X$ for some $b_{m_1}\in B_{m_1}$. Set $m_2=rm_1$.
We then have 
that 
\[ m_2 D_{m_1}= rm_1 (b_{m_1})^-_X= m_1(b^{r}_{m_1})^-_X\leq m_1 D_{m_2}.\]
This completes the proof of Lemma \ref{mono}.\\ \\
In the next section we will recall some technical tools - multivaluations and
Newton-Okounkov bodies - that wil be useful in the construction of $D({\bf B})$. 
\section{Multivaluations and Newton-Okounkov bodies.}
The proof of the various properties of the divisor $D({\bf B})$ 
will depend on the
use of multivaluations on the function field of $X({\bf B})$ and
Newton-Okounkov constructions in order to
estimate volumes of algebras using convex bodies in $\mathbb{R}^d$ associated to
admissible flags. 
\subsection{Multivaluations.}
We consider a $k$-variety $X$ of dimension $d$ equipped with an 
{\it admissible flag}.
By an admissible flag  on $X$ we mean data of a flag of varieties 
\[\hat{X}=Y_0\supset Y_1\supset Y_2\supset Y_d\] 
such that
\begin{enumerate}
\item $\hat{X}\rightarrow X$ is a birational modification. 
\item each of the varieties $Y_i$ is reduced, irreducible and of dimension $d-i$,
\item each of the varieties $Y_i$ is smooth in a neighbourhood of the point $Y_d$.
\end{enumerate}
To any such flag we can associate a multivaluation on $K(X)$, ie., a map 
\[\nu_{Y_\bullet}: K(X)\setminus 0 \rightarrow \mathbb{Z}^d\]
such that
\begin{enumerate}
\item $\nu_{Y_\bullet}(k)=0$.
\item for any $g_1, g_2 \in K(X)$ we have that $\nu_{Y_\bullet}(g_1g_2)= \nu_{Y_\bullet}(g_1) + \nu_{Y_\bullet}(g_2)$,
\item for any $g_1, g_2\in K(X)$ we have that 
\[\nu_{Y_\bullet}(g_1+g_2)\geq  {\rm min}(\nu_{Y_\bullet}(g_1) , \nu_{Y_\bullet}(g_2)),\] where the order used on 
$\mathbb{Z}^d$ is the lexicographic order. Moreover, this inequlity is an 
identity whenever $\nu_{Y_\bullet}(g_1) \neq \nu_{Y_\bullet}(g_2)$
\end{enumerate}
We recall the definition of this multivaluation.
\begin{defn}
Choose functions $f_1,\ldots, f_d\in K(X)$, 
regular in a neighbourhood of the point
$Y_d$,  such that in a neighbourhood of $Y_d$ we have for all $i$
\[ Y_i|Y_{i-1}= {\rm zero}(f_i|_{Y_{i-1}})\]
scheme-theoretically. 
Now, for any $g\in K(X)$ we define inductively functions $g_i\in K(Y_i)$ and integers $\nu_i$ by $g_0= g$ and for all $i\geq 1$ 
\begin{enumerate}
\item $\nu_i=$ degree of vanishing of $g_{i-1}$ along 
$Y_i$\footnote{which can of course be negative if $g_{i-1}$ has a pole along $Y_i$}
\item $g_i= \left(\frac{g_{i-1}}{(f_i|_{Y_{i-1}})^{\nu_{i}}}\right)|_{Y_i}$.
\end{enumerate}
The multivaluation $\nu_{Y_\bullet}$ on $K(X)$ is then defined by
\[ \nu_{Y_\bullet}(g)=(\nu_1,\nu_2,\ldots, \nu_d).\]
\end{defn}
It is immediate that this function satisfies conditions (1)- (3). 
Following the argument
in \cite{LM}, we can see that it also satisfies the condition 
\begin{lemma}\label{counting}
For any finite dimensional sub $k$-vector space $V\subset K(X)$ of finite dimension we have that
\[ \#\{ \nu_{Y_\bullet}(V\setminus 0)\}= {\rm dim}(V).\]
\end{lemma}
\noindent
{\bf Proof of Lemma \ref{counting}.} \\ \\
We proceed by induction on the dimension of $V$. The case ${\rm dim}(V)=1$
is immediate. Set $\mu= {\rm max}(\nu_{Y_\bullet}(V))$ and consider the subspace 
$V_\mu= \{ f\in V| \nu_{Y_\bullet}(f)=\mu\}$. \\ \\ 
Choose $V'$, a complement of $V_\mu$ in $V$. By the induction hypothesis, 
$\nu_{Y_\bullet}(V')={\rm dim}(V')$ and hence by (3) $\nu_{Y_\bullet}(V)=
\nu_{Y_\bullet}(V')+ \mu$. It is therefore enough to show that 
${\rm dim}(V_\mu)=1$. If not, then $V_\mu$ contains two independent 
elements $f_1$ and $f_2$ and there exists a $\lambda$ such that
\[ \nu_{Y_\bullet}(f_1+\lambda f_2)> \mu,\]
which is impossible by definition of $\mu$. This completes the proof of 
Lemma \ref{counting}.
\subsection{Construction of the Newton-Okounkov body.}
\noindent
The set of points associated by $\nu_{Y_\bullet}$ to a subspace $B$ will be an 
important invariant in what follows.
\begin{defn}
For any $k$-variety $X$, any subspace $B\subset K(X)$ and any choice of 
admissible flag $Y_\bullet$
on $X$ we denote by $\Gamma_{Y_\bullet}(B)$ the set
\[  \nu_{Y_\bullet}(B\setminus 0).\]
\end{defn}
\noindent
We note that if $B_m$ ($m\in\mathbb{N})$ is a family of subspaces of 
$K(X)$ with the property that $B_0=k$ and  $\oplus_m B_m$ is a graded algebra,
then the following set 
\[\left\{ (m, a_1,\ldots, a_d)\in \mathbb{Z}^{d+1}|\exists f\in B_m \mbox{ such that }(a_1,\ldots, a_d)= \nu_{Y_\bullet}(f)\right\}\] 
is a semigroup,  
which we denote by $\Gamma_{Y_\bullet}(\oplus_m B_m)$. We may now define the 
Newton-Okounkov body of the graded algebra $\oplus_m B_m$ with respect to the 
flag $Y_\bullet$.
\begin{defn}
Let $X$ be a $k$-variety and $(B_m)\subset K(X)$ be a family of
$k$-vector spaces such that $B_0=k$ and
$\oplus_m B_m$ is a graded algebra. Let $Y_\bullet$ be an admissible flag on $X$.
The Newton-Okounkov body of the graded algebra $\oplus_m B_m$ 
is given by
\[ \Delta_{Y_\bullet}(\oplus B_m) = \overline{{\rm Cone}(\Gamma_{Y_\bullet}(\oplus 
B_m))\cap \{(1, v)| v\in \mathbb{R}^d\}}\]
where the closure is taken with respect to the Euclidean topology on 
$\mathbb{R}^d$.
\end{defn}
\noindent
We will need the following lemma, first proved by  \cite{KK}. (The 
given here is the version given in \cite{LM}, where
proof relies heavily on previous work by Kaveh and Khovanskii). 
\begin{lemma}\label{LM}
Suppose that the valuation $\nu_{Y_\bullet}$ is such that 
the semigroup $\Gamma_{Y_\bullet}(\oplus B_m)\subset \mathbb{N}^{d+1}$ and 
moreover the following conditions are satisfied.
\begin{enumerate}
\item $\Gamma_{Y_\bullet}(B_0)=\{ 0\}$
\item there exists a finite set of vectors $(v_i,1)$ spanning a semigroup 
containing $\Gamma_{Y_\bullet}(\oplus B_m)$,  
\item $\Gamma_{Y_\bullet}(\oplus B_m)$ generates $\mathbb{Z}^{d+1}$ as a group.
\end{enumerate}
Then $\oplus_m B_m$ has a volume and moreover 
\[{\rm vol}(\oplus_m B_m)= d! {\rm vol}(\Delta_{Y_\bullet}(\oplus B_m)).\]
where the volume appearing in the right hand side is simply the standard 
Euclidean volume on $\R^d$. 
   \end{lemma}
\section{Definition of the infinite divisor $D({\bf B})$.}
We define  the infinite divisor $D({\bf B})$ 
to be the limsup of the normalisations of the divisors 
$D_m$ constructed above. 
\begin{defn}
Let ${\bf B}$ be an approximable algebra over $k$, let $X({\bf B})$ be the
associated $k$-variety and for every $m$ such that $B_m$ is 
non-trivial\footnote{ By definition of an approximable algebra, this holds for 
all large enough $m$.} let $D_m$ be the divisor defined above.
We then define  
\[D=\limsup_m\left(\frac{D_m}{m}\right).\] 
\end{defn}
We will first check that this definition makes sense, ie., that there is no prime divisor
$C\subset X$ such that $\limsup_m {\rm coeff}\left(C, 
\left(\frac{D_m}{m}\right)\right)=+\infty$. Our proof will
depend on the following proposition.
\begin{prop}\label{welldefined}
Let ${\bf B}$ be an approximable algebra over the field $k$, $X({\bf B})$
the associated algebraic variety and $Y_\bullet$ an admissible flag on $X$. 
The Newton-Okounkov body $\Delta_{Y_\bullet}({\bf B})$ is then a bounded
subset of $\R^d$.  
\end{prop}
\noindent
{\bf Proof of Proposition \ref{welldefined}.}
\\ 
We start with the following lemma.
\begin{lemma}\label{simplex} 
Let ${\bf B}$, $X({\bf B})$ and $Y_\bullet$ be as in the statement of Proposition
\ref{welldefined}. 
There is then a natural number $k$ such that the set
$\Gamma_{Y_\bullet}(B_k)$ contains a set of
$d+1$ points, $v_1,\ldots, v_{d+1}\in \R^d$, 
whose convex hull is a $d$-dimensional 
simplex of non-empty interior.  
\end{lemma}
\noindent
{\bf Proof of Lemma \ref{simplex}.}\\ \\
Choose an integer $k'$ such that $\langle B_{k'}\rangle$ is an algebra of 
positive volume, which is possible by the definition of an approximable 
algebra. 
This algebra is a sub-algebra of $\oplus_n H^0(nD_{k'})$, so by \cite{LM} 
Lemma 1.10, its Newton-Okounkov body is a compact body. Moreover,
by Lemma $(4)$, we have that $\Delta_{Y_\bullet}(\langle B_{k'} \rangle )$ has 
non-empty 
interior. The lemma follows on setting $k=rk'$ for sufficiently large $r$.
This completes the proof of Lemma \ref{simplex}.
\\ \\
We return now to the proof of Proposition \ref{welldefined}. 
For any $k$ we have that 
\[ \Delta_{Y_\bullet}({\bf B}_k)= k \Delta_{Y_\bullet}({\bf B})\]
so after passing to the algebra ${\bf B}_{k}$ we may assume that 
$k=1$. We denote by $v_1,\ldots, v_{d+1}$ the vectors whose existence is 
guaranteed by Lemma \ref{simplex}. \\ \\
Note that for any $m$ the map
\[ \left\{(a_1,\ldots, a_{d+1})| a_i\geq 0, \sum_{i=1}^{d+1} a_i =m\right\}\rightarrow  \sum_{i=1}^{d+1} a_i v_i\]
is injective. We set
\[ M=\sup_i \{||v_i ||\}\]
Suppose that $\Delta_{Y_\bullet}({\bf B})$ is not bounded. In other words, 
for any $N\in \R^+$ there exists a natural number $m(N)$ and 
an element
$b\in B_{m(N)}$ such that $|| \nu(b)|| > N m(N)$. 
For simplicity we will consider only values of $N$ which are exact 
multiples of $M$. \\ \\
Let us now consider for any $p$
the set $\Gamma_{Y_\bullet}(B_{pm(N)})$: our aim is to show that this set is too large.
Note that  $\Gamma_{Y_\bullet}(B_{pm(N)})$ contains the set
\[ Z(p, m(N))= \left\{ a_0\nu(b)+\sum_{i=1}^{d+1}a_i v_i| a_i\geq 0, m(N)a_0+ \sum_{i=1}^{d+1} a_i= pm(N)\right\}.\]
Once again, for the sake of simplicity we assume that $p$ is a multiple of $N$. The map  
\[ (a_0,\ldots, a_{d+1}) \rightarrow  a_0\nu(b)+\sum_{i=1}^{d+1}a_i v_i\]
is not a priori injective, but it becomes injective if we make the additional 
assumption that 
that $a_0$ be a multiple of $2pM/N$. It follows that
\[ {\rm dim}(B_{pm(N)})= \#\Gamma_{Y_\bullet}(B_{pm(N)})\]
\[\geq \#\left\{ (a_0,\ldots, a_{d+1})| a_i\geq 0, a_0\mbox{ is a multiple of 2pM/N },\;m(N)a_0+ \sum_{i=1}^{d+1} a_i= pm(N)\right\}.\]
To find a lower bound on  ${\rm dim}(B_{pm(N)})$ it will be enough to find a lower bound on the size of the set
\[ \left\{ (a_0,\ldots, a_{d+1})| a_i\geq 0, a_0\mbox{ is a multiple of 2pM/N },
\;m(N)a_0+ \sum_{i=1}^{d+1} a_i= pm(N)\right\}.\]
On setting $c= \frac{a_0 N}{2pM}$ we see that the size of this set is equal to
\[ \sum_{c=0}^{N/2M} {pm-2pmMc/N +d \choose d}
\geq \sum_{c=0}^{N/2M}\frac{1}{d!}(p^dm^d(1-2Mc/N)^d)\]
\[= \frac{1}{d!} (p^d m^d)\sum_{c=0}^{N/2M} (1-2Mc/N)^d 
\geq \frac{1}{d!} (p^d m^d)\sum_{c=0}^{N/4M} (1-2Mc/N)^d\]
\[ \geq \frac{Np^dm^d}{2^{d+2}d!M}\]  
from which it follows that 
\[ {\rm vol}({\bf B})\geq \frac{N}{M 2^{d+2}}.\]
 for any $N$. This is impossible since ${\bf B}$ has finite volume.
This completes the proof of Proposition \ref{welldefined}.\\ \\
Given this proposition, it is fairly easy to deduce the following.
\begin{prop}\label{wd}
Let ${\bf B}$ and $X({\bf B})$ be as above. The divisor $D({\bf B})$ 
constructed above is then well-defined, ie., for any prime divisor $C\subset X$
the set
\[ \{ \mbox{ coeff }(C; D_m/m)| m\in \mathbb{N}\}\]
is bounded.
\end{prop}
\noindent
{\bf Proof of Proposition \ref{wd}.}
\\ \\
We argue by contradiction : suppose that $C$ is a divisor on $X$ such that 
 $\mbox{ coeff }(C; D_m/m)$ can be arbitrarily large. 
After blow-up, we may assume that $C$ is a smooth divisor in $X$ and 
choose an admissible flag $Y_\bullet$ on $X({\bf B})$ whose first element is $C$ - the 
Newton-Okounkov body of ${\bf B}$ with respect to $Y_\bullet$ is then unbounded,
contradicting Proposition \ref{welldefined}. 
This completes the proof of Proposition \ref{wd}. \\ \\
In order to have a reasonable level of control of the algebra $R(D({\bf B}))$
it will be useful to have the following result.
\begin{prop}\label{finite}
Let ${\bf B}$, $X({\bf B})$ and $D({\bf B})$ be as above.
For any $m\in \mathbb{N}$ the divisor $\lfloor mD\rfloor$ is in fact 
a finite divisor. 
\end{prop}
\noindent
In other words, if $D({\bf B})=\sum_i a_i D_i$ then 
the sequence $a_i$ converges to $0$.
\\ \\ 
{\bf Proof of Proposition \ref{finite}.}\\ \\
We know that for $k$ sufficiently large the map associated to the
linear system $B_k$ is birational. Note that if $\lfloor mD \rfloor$ is not 
finite then $\lfloor mkD \rfloor$ is not finite either, so on passing to 
the approximable algebra ${\bf B}_k$ we may assume that $k=1$. \\
Let $U$ be the open set on which this birational map associated to the linear
system $B_1$ is an isomorphism. 
Proposition \ref{finite} will follow from the following result. 
\begin{prop}\label{bound}
Assume that ${\bf B}$ is an approximable algebra such that $B_1$ gives rise to
 a birational map on $X({\bf B})$.  Let $U$ be the set on which this map 
is an isomorphism. 
Choose $m$ such that \[{\rm vol}(\langle B_m\rangle )+
\frac{ \delta^{d+1}}{(1+\delta)^d}\geq 
{\rm vol}({\bf B}).\] Let $C$ be a prime divisor that does not appear in 
$D_{m'}$ for any $m'\leq m$ and which is not contained in $X\setminus U$.
Then for all integers $m''$ the coefficient of $C$ in $D_{m''}$ is less than
$\delta m''$.   
\end{prop} 
\noindent
{\bf Proof of Proposition \ref{bound}.}
\\ \\
We argue by contradiction. Suppose that the conclusion of the proposition is 
false and that there exists an $m''$ such that the coefficient of $C$ in 
$D_{m''}$ is $>\delta m''$. Since it then follows that the coefficient of $C$ in $D_{rm''}$ is $> \delta r m''$ for all $r$, we may assume that $m''$ is a 
multiple of $m$.  

Choose a point $x$ contained in $C$ which satisfies the two following criteria.
\begin{enumerate}
\item $x$ is contained in $U$
\item $x$ is not contained in $D_{m'}$ for any $m'\leq m$
\end{enumerate}
It is possible to choose such an $x$ because of the assumptions on $C$. Now,
since $X$ is defined up to birational equivalence, we may blow up the point 
$X$ and
consider a generic 
infinitesimal flag centred at the point $x$, ie. a flag $(Y_1,\ldots, Y_d)$
such that $Y_1$ is the exceptional divisor of $X$ over $x$ and
the $Y_i$s for $i=2\ldots d$ are very general linear subspaces of $Y_1$. \\ \\
The various conditions required on the point $x$ have the following 
implications.
\begin{enumerate}
\item $\Gamma_{Y_\bullet}({\rm B_1})$ contains the vectors 
$e_1=(0,0,\ldots, 0)$, $e_2=(1,0\ldots, 0)$, $e_3=(1,1,\ldots, 0)$ , 
$(1,0,1,\ldots,0)$, $e_{d+1}=(1,0,\ldots,0, 1)$ because the linear system $B_1$
defines an isomorphism in a neighbourhood of $x$. 
\item $\Gamma_{Y_\bullet}(\langle B_m\rangle )$ is contained in 
$\mathbb{N}^{d+1}$
because $x\not\in D_m$.
\item If $C$ appears in $D_{m''}$ with coefficient $\geq \delta m''$ then 
$\Gamma_{Y_\bullet}(B_{m''})$ contains a vector $v=(v_1,0,\ldots, 0)$ where 
$v_1\leq -\delta m''$. (The subsequent coefficients of $v$ are all 
zero because of the very general choice of the linear subspaces $Y_i$). 
\end{enumerate}
Let $r$ be the integer such that $m''=rm$. 
Now, let us consider for large enough $p$ the set
$\Gamma_{Y_\bullet}( B_{pm''})$, and compare it with the set
$\Gamma_{Y_\bullet}(S^{rp}B_m)$. Recall that these two groups are supposed to 
be about the same size.
\\ \\
The set $\Gamma_{Y_\bullet}(B_{pm''})$ contains the set of vectors 
\[ \left\{a_0 v+ \sum_{i=1}^{d+1} a_i e_i| a_i\geq 0,  am''+ 
\sum_{i=1}^{d+1} a_i =pm''.\right\}\]
We try to bound below the number of elements of this set that are
not contained in  $\Gamma_{Y_\bullet} (S^{rp}B_m)$.
We know
that $\Gamma_{Y_\bullet} ({S}^{rp}B_m)$ does not contain any vectors with
negative first coefficient. 
It follows that the set \[\Gamma_{Y_\bullet}(B_{pm''})\setminus 
\Gamma_{Y_\bullet}((S^{rp}B_m))\] 
contains the following set of integral vectors
\[ \left\{ a_0 v+ \sum_{i=1}^{d+1} a_i e_i|a_i\geq 0,  a_0 \delta m''
\geq \sum_{i=2}^{d+1} a_i,a_0m''+ \sum_{i=1}^{d+1} a_i =pm'' \right\}.\]
In presence of the equation $a_0m''+\sum_{i=1}^{d+1} a_i =pm''$ the inequality 
$ a_0 \delta m''\geq \sum_{i=2}^{d+1} a_i$ is satisfied
whenever 
\[ p\geq a_0 \geq \frac{p}{\delta+1}.\]
So to obtain a lower bound on 
${\rm vol}({\bf B})-{\rm vol} \langle B_{m''}\rangle $ it will be enough to
obtain a lower bound on the size of the set
\[S(m'',p, \delta)
=  \left\{ a_0 v+ \sum_{i=1}^{d+1} a_i e_i|a_i\geq 0, a_0m''+ \sum_{i=1}^{d+1} 
a_i =pm'', p\geq a_0 \geq \frac{p}{\delta+1} \right\}.\]
Let us study first the tuples $(a_0,\ldots, a_{d+1})$ and $(a'_0,
\ldots, a'_{d+1})$ such that 
\[ a_0 v+ \sum_{i=1}^{d+1} a_i e_i=a'_0 v+ \sum_{i=1}^{d+1} a'_i e_i\]
satisfying the condition 
\[a_0m''+ \sum_{i=1}^{d+1} a_i =a'_0m''+ \sum_{i=1}^{d+1} a'_i =pm''\]
On looking at the $i$ th  coordinate of this vector we see that $a_i= a_i'$
for $i\geq 3$. We therefore have that
\[ a_0v_1+a_2= a_0'v_1+a'_2\]
so that $a_2-a_2'$ is a multiple of $|v_1|> \delta m''$. In particular, if 
we impose the extra condition $a_2, a_2'\leq\delta m''$ then 
$a_2=a_2'$.
Note that the condition $a_0m''+ \sum_{i=1}^{d+1} a_i =a'_0m''+ \sum_{i=1}^{d+1} 
a'_i $ implies that if 
$a_0=a'_0$ and $a_i=b_i'$ for $i\geq 2$ then $a_1=a'_1$. \\ \\
Pulling the above together we see that on the set
\[ Z(m'',p,\delta)=\left\{(a_0,\ldots, a_{d+1})|   a_i\geq 0, a_0m''+ 
\sum_{i=1}^{d+1} a_i =pm'', p\geq a_0 \geq \frac{p}{\delta+1}, 
0 \leq a_2< \delta m''\right\}\]
the map
$(a_0,\ldots, a_{d+1})
\rightarrow  a_0 v+ \sum_{i=1}^{d+1} a_i e_i$
is injective, so the size of $\Gamma_{Y_\bullet}(B_{pm''})\setminus 
\Gamma_{Y_\bullet}((S^{rp}B_m))$ is bounded below by the size of $Z(m'',p,\delta)$. 
Now, we have that for large enough $p$ 
\[ \# Z(m'',p,\delta)= \sum_{a_0=\frac{p}{\delta +1}}^p 
\sum_{a_2=0}^{\delta m''} { (pm''-a_0m''-a_2) +d-1\choose d-1}\]
\[ \geq \sum_{a_0=\frac{p}{\delta+1}}^{p-\delta} \sum_{a_2=0}^{\delta m''} \frac{1}{(d-1)!} (((pm''-am''-a_2))^{d-1}\]
\[ \geq \sum_{a_0=\frac{p}{\delta +1}}^{p-\delta} \sum_{a_2=0}^{\delta m''} \frac{1}{(d-1)!} (((pm''-am''-\delta m''))^{d-1}\]
\[=\sum_{a_0=\frac{p}{\delta +1}}^{p-\delta} \delta m'' \frac{1}{(d-1)!} (((pm''-am''-\delta m'')^{d-1}\]
\[= \sum_{a'=\delta}^{\frac{p\delta}{\delta +1}} \frac{\delta m''}{(d-1)!} 
((a'-\delta)m'')^{d-1}
=\frac{\delta  m''^d  }{(d-1)!}\sum_{a'=\delta}^{\frac{p\delta }{\delta +1}}(a'-\delta )^{d-1}\]
\[ \sim \frac{ \delta m''^d}{(d-1)!}\frac{1}{d} \left(\frac{p\delta }{\delta+1}-\delta \right)^d 
\sim \frac{ \delta^{d+1}  (pm'')^d}{d! (1+\delta)^d}.\]
In particular, 
\[ {\rm dim} (B_{pm''})-{\rm dim} (S^{rp}(B_m))\geq  \frac{ \delta^{d+1} (pm'')^d}{d! (1+\delta )^d}\]
\[ \limsup_{p\rightarrow \infty} \frac{ {\rm dim} (B_{pm''})-{\rm dim} ({\rm S}^{rp}(B_m))}{(pm'')^d/d!}\geq  \frac{ \delta^{d+1} }{(1+\delta )^d}\]
\[ {\rm vol}({\bf B}) -{\rm vol}(\langle
B_m\rangle)\geq \frac{ \delta^{d+1}}{d! (1+\delta )^d}.\]
But this contradicts the initial assumptions on $m$. 
This completes the proof of Proposition \ref{bound}.\\ \\
\noindent The proof of Proposition \ref{finite} follows easily. Choose an 
integer $l$ and set $\delta=1/l$. Let $m$ be such that  ${\rm vol}(\langle B_m\rangle )+
\frac{ \delta^{d+1}}{(1+\delta)^d}\geq 
{\rm vol}({\bf B}).$ By Proposition \ref{bound},
any divisor $C$ which appears in 
$\lfloor lD\rfloor$, is either contained in $X\setminus U$ or in 
$\cup_{m'\leq m} D_{m'}$. There is only a finite number of such 
divisors, so this completes the proof of Propositon \ref{finite}.\\ \\
\noindent
Pulling the above together, we have the following theorem.
\begin{thm}\label{mainthm}
Let ${\bf B}=\oplus_m B_m$ be an approximable algebra over a field $k$
of dimension $d({\bf B})$. Then there exists a smooth
variety $X({\bf B})$ of dimension $d({\bf B})$ and an infinite 
Weil divisor on 
$X({\bf B})$ which we denote $D=\sum_i a_i D_i$ such that 
\begin{enumerate}
\item for all $m$ the divisor $\lfloor mD \rfloor$ is a finite sum of prime
divisors, 
\item there is an inclusion of graded algebras
\[ {\bf B}\hookrightarrow \oplus_m H^0(X, \mathcal{O}_X(\lfloor mD \rfloor))\]
\end{enumerate}
\end{thm}
\noindent
In the final section, we will prove a partial result in the opposite direction.
\section{Approximable subalgebras of  $\oplus_m H^0(X, \mathcal{O}_X(\lfloor mD \rfloor))$.}
\noindent
Given \ref{mainthm} it seems reasonable to wonder when a subalgebra of 
$\oplus_m H^0(X, \lfloor mD \rfloor)$, where $D$ is
an infinite sum of Weil divisors on a variety $X$ over an algebraically closed field of characteristic zero, is an approximable algebra. 
We start with the following result
\begin{prop}\label{approx}
Let ${\bf C}=\oplus_m C_m$ be an approximable algebra of dimension $d({\bf C})$
over an algebraically closed field $k$ of characteristic zero.
Let ${\bf B}=\oplus_m B_m$ be a sub graded algebra of ${\bf C}$ such that 
\[ \left( \frac{{\rm rk}(B_m)}{ m^{d({\bf C})}}\right)\]
does not converge to $0$ and
for all sufficiently large $n$ the space $B_n$ is non-empty. 
Then the algebra ${\bf B}$ is also approximable.
\end{prop}
Throughout the proof of this proposition we will denote$d({\bf C})$ by $d$. \\ \\
{\bf Proof of Proposition \ref{approx}.} 
\\ \\ Conditions 1) and 2) of the 
definition of an approximable algebra are immediately satisfied, so it remains
to check only the condition (3). 
Pick an admissible flag $Y_\bullet$ on $X=X({\bf C})$
such that $Y_d$ is not contained in $D({\bf C})$, so that $\Gamma_{Y_\bullet}(C_m)\subset \mathbb{N}^d$ for all $m$. 
We know that the Newton-Okounkov 
body of ${\bf C}$ is compact. It follows that the Newton-Okounkov body of 
${\bf B}$ is also compact and since  
$\left( \frac{{\rm rk}(B_m)}{ m^{d({\bf C})}}\right)$ does not converge to $0$ 
the Newton-Okounkov body of ${\bf B}$ has non-empty interior. \\ \\
There is therefore a $k$ such that $B_k$ contains $(d+1)$ elements $(f_0,
\ldots, f_d)$ whose valuation vectors span a $(d+1)$-dimensional simplex in 
$\R^d$. It follows that the elements $\frac{f_1}{f_0}\ldots, \frac{f_d}{f_0}$
are algebraically independent in $K^{\rm hom}({\bf C})$ so the induced rational 
map $X\rightarrow \mathbb{P}^d$
has dense image.\\ \\
Let $x$ be a point of $X$ at which this rational morphism is well-defined and 
immersive, and which is not contained in $D_i$ for any $i$.
Let $Y'_\bullet$ be an infinitesimal flag centered at the point $x$. 
We then know that the set $\Gamma_{Y_\bullet}(B_{rk})$ generates $\mathbb{Z}^d$
as group for any $r$ and it follows that for any large enough $n$ the semigroup
$\Gamma_{Y\bullet}\langle B_n \rangle$ generates $\mathbb{Z}^{d+1}$. The conditions of
Lemma \ref{LM} therefore apply so we can say that ${\rm vol}({\bf B})$ exists, as does ${\rm vol}(\langle B_m\rangle)$ for any large enough $n$ and moreover
\[ {\rm vol}(\langle B_n \rangle)= d! {\rm vol}(\Delta_{Y_\bullet}(\langle B_n \rangle)).\]
Equally, we have that 
\[ {\rm vol}({\bf B})=  d!{\rm vol}(\Delta_{Y_\bullet}({\bf B})).\]
We have only therefore to prove that \[\lim_n {\rm vol}(\Delta_{Y_\bullet}(\langle B_n \rangle))/n^d ={\rm vol} (\Delta_{Y_\bullet}({\bf B})).\]
It is immediate that \[\limsup_n {\rm vol}(\Delta_{Y_\bullet}(\langle B_n \rangle))/n^d ={\rm vol} (\Delta_{Y_\bullet}({\bf B}))\] 
so it remains only to prove that this $\limsup$ is in fact a limit. But we know that for $k>k_0$ 
\[ {\rm vol}((\Delta_{Y_\bullet}(\langle B_{rn+k} \rangle))/(rn+k)^d\geq  \left(\frac{rn}{rn+k}\right)^d{\rm vol}((\Delta_{Y_\bullet}(\langle B_{rn}\rangle))/(rn)^d\]
\[\geq \left(\frac{rn}{rn+k}\right)^d
{\rm vol}((\Delta_{Y_\bullet}(\langle B_{n}\rangle))/(n)^d\]
so for any $n$ we have that
\[ \liminf_m {\rm vol}((\Delta_{Y_\bullet}(\langle B_{m}\rangle))/m^d\geq  {\rm vol}((\Delta_{Y_\bullet}(\langle B_{n}\rangle))/(n)^d\]
so that $\liminf_m {\rm vol}((\Delta_{Y_\bullet}(\langle B_{m}\rangle))/m^d= \limsup_m {\rm vol}((\Delta_{Y_\bullet}(\langle B_{m}\rangle))/m^d$ and 
\[ \lim_n {\rm vol}(\Delta_{Y_\bullet}(\langle B_n \rangle))/n^d ={\rm vol} (\Delta_{Y_\bullet}({\bf B})).\]
This is completes the proof of Proposition \ref{approx}. 
\noindent
\begin{prop}\label{approx2}
Let $X$ be a smooth algebraic variety defined on an algebraically closed field of charcteristic zero.  
If $D=\sum a_iD_i$ is an infinite Weil divisor on $X$ such that the class $[D]= \sum_i a_i[D_i]$ converges towards a finite numerical big divisor class in $NS(X)$ 
then the ring ${\bf B}=\oplus_mH^0(\lfloor mD \rfloor)$ is approximable.
\end{prop}
{\bf Proof of Proposition \ref{approx2}.}
Set $D_m= \lfloor mD\rfloor$ and 
$B_m= H^0(\lfloor mD \rfloor/m)$, so that ${\bf B}=\oplus_m B_m$. 
Conditions $(1)$ and $(2)$ of the definition of an approximable algebra are 
immediately satisfied so our aim is to show that 
condition $(3)$ is also satisfied.\\ \\ 
We have that $K^{\rm hom}({\bf B})=K(X)$ since for large enough $m$ the divisor 
$D_m$ is big. After choosing an element $f\in B_1$ we have an inclusion
\[ K^{\rm hom}({\bf B})\hookrightarrow K(X).\]
We pick $Y_\bullet$, a very general infinitesimal flag on $X$ centred at 
$y\in X$, which for this reason is not contained in $D_i$ for any $i$. We can 
therefore calculate the Newton-Okounkov body of ${\bf B}$ with respect to
$Y_\bullet$. 
Lemma 1.10 of Lazarsfeld and Mustata shows that $\Delta_{Y_\bullet}({\bf B})$ is 
compact.\\ \\
Let us show further that ${\bf B}$ satisfies the conditions (1)-(3)
 of Lemma 5. Condition (1) holds by definition and condition (2) is immediate 
because the Newton Okounkov body is bounded: it remains only to prove condition
(3). 
We note that
$(3)$ holds for ${\bf B}$ if it holds for $\langle B_m\rangle$ for some 
$m$. \\ \\ 
Since the class $[D]$ is big, the class $\lfloor mD \rfloor$ is big for 
large enough $m$ so for large enough $m$ the map defined by $B_m$
is birational, from which it follows that condition $(3)$ is satisfied for
${\bf B}$. It follows that ${\bf B}$ has a well-defined volume,
equal to the volume of its Newton-Okounkov body with respect to a general 
infinitesimal flag.
Moreover, for all large enough $m$ 
\[{\rm vol}(\langle B_m\rangle)= d! {\rm vol}(\Delta_{Y_\bullet}(\langle B_m\rangle)).\] 
It will be enough to show that
\[ \lim_m {\rm vol}(\Delta_{Y_\bullet}(B_m)/m^d)= {\rm vol}({\bf B}).\]
It is immediate that 
\[ \limsup_m {\rm vol}(\Delta_{Y_\bullet}(\langle B_m\rangle))/m^d= {\rm vol}({\bf B}).\]
so it remains only to prove that this $\limsup$ is in fact a limit. But we know that for $k>k_0$ 
\[ {\rm vol}((\Delta_{Y_\bullet}(\langle B_{rn+k} \rangle))/(rn+k)^d\geq  \left(\frac{rn}{rn+k}\right)^d{\rm vol}((\Delta_{Y_\bullet}(\langle B_{rn}\rangle))/(rn)^d\]
\[\geq \left(\frac{rn}{rn+k}\right)^d
{\rm vol}((\Delta_{Y_\bullet}(\langle B_{n}\rangle))/(n)^d\]
so for any $n$ we have that
\[ \liminf_m {\rm vol}((\Delta_{Y_\bullet}(\langle B_{m}\rangle))/m^d\geq  {\rm vol}((\Delta_{Y_\bullet}(\langle B_{n}\rangle))/(n)^d\]
so that $\liminf_m {\rm vol}((\Delta_{Y_\bullet}(\langle B_{m}\rangle))/m^d= \limsup_m {\rm vol}((\Delta_{Y_\bullet}(\langle B_{m}\rangle))/m^d$ and 
\[ \lim_n {\rm vol}(\Delta_{Y_\bullet}(\langle B_n \rangle))/n^d ={\rm vol} (\Delta_{Y_\bullet}({\bf B})).\]
This completes the proof of the proposition.\\ \\
Pulling these two results together, we see that we have the following.
\begin{thm}
Let $X$ be a complex algebraic variety and let $D=\sum a_i D_i$ be an infinite 
Weil divisor on $X$ such that the sum of divisor classes 
$ \sum_i a_i[D_i ]$ converges
to a finite real big cohomology class. Any graded subalgebra of 
$\oplus_m H^0(mD)$ such that 
\[ \left( \frac{{\rm rk}(B_m)}{ m^{d({\bf C})}}\right)\]
does not converge to zero is then an approximable algebra.
\end{thm}
\begin{rem}
The following intriguing question remains open : do there exist approximable
algebras ${\bf B}$ such that the associated infinite Weil 
divisor $D({\bf B})=\sum_i a_i D_i$ 
does not converge in the space of numerical classes of divisors ?
\end{rem}

Author - Catriona Maclean\\ \\
Adress - Institut Fourier, 100 rue des maths, Saint Martin d'H\`eres, 38400, France \\ \\
Email - Catriona.Maclean@univ-grenoble-alpes.fr\\ \\
Affiliation - Universit\'e Grenoble Alpes.
\end{document}